\input amstex
\documentstyle{amsppt}
\pagewidth{4.50truein}
\pageheight{7.75truein}
\hoffset=1.0truein
\voffset=0.125truein

\input xy
\xyoption{matrix}\xyoption{arrow}\xyoption{curve}\xyoption{frame}

\def\edge{\ar@{-}}

\long\def\ignore#1{#1}
\long\def\ignorethree#1{#1}

\def\pinf{{\Cal P}^\infty}
\def\lamod{\Lambda{\operatorname{-mod}}}
\def\laMod{\Lambda{\operatorname{-Mod}}}
\def\dlmod{\Delta{\operatorname{-mod}}}
\def\rmod{R{\operatorname{-mod}}}
\def\lfindim{\operatorname{l\,fin\,dim}}
\def\lFindim{\operatorname{l\,Fin\,dim}}
\def\A{{\Cal A}}
\def\B{{\Cal B}}
\def\C{{\Cal C}}
\def\arrA{\overarrow{\Cal A}}
\def\arrpinf{\overarrow{\pinf}}
\def\Hom{\operatorname{Hom}}
\def\kgam{{K\Gamma}}
\def\Soc{\operatorname{Soc}}
\def\pdim{\operatorname{p\,dim}}
\def\tilgam{\widetilde{\gamma}}
\def\tilbet{\widetilde{\beta}}
\def\NN{{\Bbb N}}
\def\bfk{{\bold k}}
\def\bfl{{\bold l}}
\def\bfm{{\bold m}}
\def\tlg{\widetilde g}
\def\abar{\overline a}
\def\bbar{\overline b}
\def\cbar{\overline c}

\topmatter

\title Co-- Versus Contravariant Finiteness of Categories of Representations
\endtitle

\rightheadtext{Co-- Versus Contravariant Finiteness}

\author B. Huisgen--Zimmermann and S. O. Smal\o\endauthor

\address Department of Mathematics, University of California, Santa Barbara,
CA 93106, USA\endaddress
\email birge\@math.ucsb.edu\endemail

\address Department of Mathematics, The Norwegian University for Science and
Technology, 7055 Dragvoll, Norway\endaddress

\email sverre\@matstat.unit.no\endemail

\thanks The first author was partially supported by a grant from the NSF
while this work was done, while the second author was supported by the
US--Norway Fulbright Foundation. Most of this research was carried out
while the second author was visiting the University of California at Santa
Barbara, and he wishes to thank his coauthor for her kind hospitality.
\endthanks

\abstract This article supplements recent work of the authors. (1) A criterion for failure of covariant finiteness of a full subcategory of $\Lambda\text{-mod}$ is given, where $\Lambda$ is a finite dimensional algebra. The criterion is applied to the category ${\Cal P}^{\infty}(\Lambda\text{-mod})$ of all finitely generated $\Lambda$-modules of finite projective dimension, yielding a negative answer to the question whether ${\Cal P}^{\infty}(\Lambda\text{-mod})$ is always covariantly finite in $\Lambda\text{-mod}$. Part (2) concerns contravariant finiteness of ${\Cal P}^{\infty}(\Lambda\text{-mod})$. An example is given where this condition fails, the failure being, however, curable via a sequence of one-point extensions. In particular, this example demonstrates that curing failure of contravariant finiteness of ${\Cal P}^{\infty}(\Lambda\text{-mod})$ usually involves a tradeoff with respect to other desirable qualities of the algebra.
\endabstract

\endtopmatter

\document

\head 1. Introduction and Terminology\endhead

Functorial finiteness conditions for certain categories of finitely
generated representations of an algebra may have a major impact also on
the non-finitely generated representations, as was shown by the authors in
\cite{10}. More precisely: Let $\Lambda$ be an artin algebra, and let
$\pinf(\lamod)$ and $\pinf(\laMod)$ be the full subcategories of the
categories $\lamod$ of finitely generated left $\Lambda$-modules and the
full module category $\laMod$, respectively, consisting of the objects of
finite projective dimension in either case. Then contravariant finiteness of
$\pinf(\lamod)$ in $\lamod$ forces arbitrary modules in $\pinf(\laMod)$ to be
direct limits of objects in $\pinf(\lamod)$. When combined with a theorem of
Auslander and Reiten \cite{1}, this entails that $\lfindim \Lambda=
\lFindim \Lambda= \sup _{1\le i \le n} \pdim A_i$ in this case, where
$A_1,\dots,A_n$ are the minimal right $\pinf(\lamod)$-ap\-prox\-i\-ma\-tions
of the simple left
$\Lambda$-modules. (Here $\lfindim \Lambda$ and $\lFindim \Lambda$ are the
suprema of the projective dimensions attained on $\pinf(\lamod)$ and
$\pinf(\laMod)$, respectively.)

As a byproduct of the described connections, one obtains that
\underbar{contra}\-var\-i\-ant finiteness of $\pinf(\lamod)$ implies
\underbar{co}variant finiteness of this category in $\lamod$. Indeed, by
Crawley-Boevey's \cite{4, Theorem 4.2}, an additive subcategory $\A$ of
$\lamod$ is covariantly finite if and only if the closure $\arrA$ of $\A$
under direct limits is closed under direct products as well. As
explained above, contravariant finiteness of $\pinf(\lamod)$ implies
$\pinf(\laMod)= \arrpinf(\lamod)$ and $\lFindim \Lambda <\infty$, which
 guarantees that $\arrpinf(\lamod)$ is closed under direct products.

It is not hard to find examples demonstrating that, in general,
contravariant finiteness of $\pinf(\lamod)$ in $\lamod$ is properly stronger
than covariant finiteness. In fact, the initial example -- due to Igusa,
Smal\o, and Todorov \cite{11} -- of a situation where $\pinf(\lamod)$ fails
to be contravariantly finite already serves to show this. This leaves one
wondering whether $\pinf(\lamod)$ might always be covariantly finite in
$\lamod$. The answer is negative, as we show here, but examples are somewhat
harder to come by.

The first part of the present paper is devoted to developing a criterion for
failure of covariant finiteness of $\pinf(\lamod)$ and to then applying it
to a finite dimensional special biserial algebra (for a  definition see
under `Notation and Terminology' below). We believe that, in the restricted
setting of special biserial algebras $\Lambda$, the conditions of this
criterion actually provide an equivalent description for failure of covariant
finiteness of
$\pinf(\lamod)$. The criterion can be considered as a somewhat weaker twin
sibling of the sufficient condition for failure of contravariant finiteness
developed by Happel and the first author in \cite{7, Criterion 10}. We
remark that, while the concepts of contra- and covariant finiteness are
mutually dual, the theories relating them to a prescribed subcategory of
$\lamod$ are of course not; in particular, the argument backing up the
criterion presented here differs substantially from that used to prove
\cite{7, Criterion 10}.

The homological picture available for algebras $\Lambda$ having the property
that $\pinf(\lamod)$ is contravariantly finite in $\lamod$ naturally raises
the question as to how abundant they are. While this condition is known to
`slice diagonally' through the standard classes of algebras, the authors
conjecture that failure can `often' be fixed in the following sense: Namely,
that there exists a sequence $\Lambda_0= \Lambda, \Lambda_1, \dots,
\Lambda_m$ of Artin algebras such that each $\Lambda_i$ is a one-point
extension of $\Lambda_{i-1}$ and $\pinf(\Lambda_m\text{-mod})$ is
contravariantly finite in $\Lambda_m\text{-mod}$. Since in each passage from
$\Lambda_{i-1}$ to $\Lambda_i$ both the little and big finitistic dimensions
increase by at most 1, this might give a handle on the difference $\lFindim
-\lfindim$ for special classes of algebras. Moreover, the existence of a
sequence as above guarantees that $\lFindim \Lambda <\infty$. In general, this
process will force one to leave a given `nice' class of algebras, however.
The second part of this article is devoted to the explicit construction of
such a sequence
$\Lambda_0= \Lambda,
\Lambda_1,
\dots,
\Lambda_m$ such that $\Lambda$ is a monomial relation algebra while
$\Lambda_m$ cannot be chosen from this class of algebras.

\subhead Terminology and Notation\endsubhead

In the following, $\Lambda$ will be a split finite dimensional algebra over
a field $K$, i.e., $\Lambda$ will be of the form $\kgam/I$ for some quiver
$\Gamma$ and an admissible ideal $I$ of the path algebra $\kgam$. The
Jacobson radical of $\Lambda$ will be denoted by $J$. For us, the primitive
idempotents of $\Lambda$ will be those which naturally correspond to the
vertices of $\Gamma$; in fact, we will identify a complete set of primitive
idempotents of
$\Lambda$ with the vertices of the quiver. 

According to \cite{2}, a full subcategory $\A$ of $\lamod$ is said to be
{\it contravariantly finite} in $\lamod$ if, for each module $M$ in
$\lamod$, there exists a homomorphism $f : A\rightarrow M$ with $A\in\A$
such that the following sequence of contravariant functors is exact:
$$\Hom_\Lambda(-,A)|_\A @>\Hom_\Lambda(-,f)>> \Hom_\Lambda(-,M)|_\A @>>> 0;$$
in other words, exactness of this sequence means that each  $g\in
\Hom_\Lambda(B,M)$ with $B\in\A$ factors through $f$. In that case, $A$ is
called a {\it right $\A$-ap\-prox\-i\-ma\-tion} of $M$. By \cite{2}, the
$\A$-ap\-prox\-i\-ma\-tions of $M$ that have minimal dimension are pairwise
isomorphic. This justifies reference to {\it the} minimal right
$\A$-ap\-prox\-i\-ma\-tion of $M$, whenever $\A$ is contravariantly finite.
The terms {\it covariantly finite} and {\it left
$\A$-ap\-prox\-i\-ma\-tion} are defined dually.

An algebra $\Lambda$ is said to be {\it special biserial} in case it is
isomorphic to a path algebra modulo relations, $\kgam/I$, with the following
properties: Given any vertex
$e$ of
$\Gamma$, at most two arrows enter $e$ and at most two arrows leave $e$;
moreover, for any arrow
$\alpha$ in $\Gamma$ there is at most one arrow $\beta$ with
$\alpha\beta\notin I$ and at most one arrow $\gamma$ with $\gamma\alpha
\notin I$.

Recall, moreover, that given two algebras $\Lambda= \kgam/I$ and $\Lambda'=
K\Gamma'/I'$, the second is called a {\it one-point extension} of the first
in case $\Gamma'$ results from $\Gamma$ through addition of a single vertex
which is a {\it source} of $\Gamma'$ such that  $I'\cap \kgam =I$. For
importance and properties of one-point extensions, we refer to \cite{13}.

Given paths $p$ and $q$ of $\Gamma$, we say that $q$ is a {\it subpath} of
$p$ if $p= p_2qp_1$ in $\kgam$ for suitable paths $p_1$ and $p_2$. We call
$q$ a {\it right (left) subpath} of $p$ in case $p= p_2q$ 
(respectively, $p= qp_1$) for suitable paths $p_2$ (respectively, $p_1$). The
path
$p$ is said to {\it start (end)} in the arrow
$\alpha$ if $\alpha$ is a right (left) subpath of $p$. Furthermore, we call
an element $x$ of a left $\Lambda$-module $M$ a {\it top element} in case
$x\in M\setminus JM$ and $ex=x$ for a primitive idempotent $e$ from our
distinguished set.

Finally, we refer the reader to previous work of the authors for their
graphing conventions (see, e.g., \cite{7,8,9,10}). The graphs most
crucial to the present note are zigzags of the type

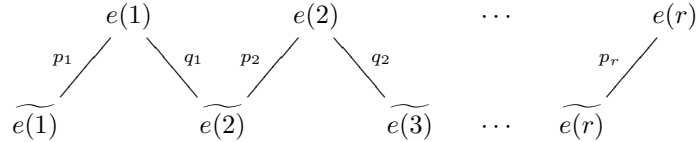
\ignore{
$$\xymatrixcolsep{1pc}
\xy\xymatrix{
 &e(1) \edge[dl]_{p_1} \edge[dr]^{q_1} &&e(2) \edge[dl]_{p_2}
\edge[dr]^{q_2} &&\cdots &&e(r) \edge[dl]_{p_r}\\
{\widetilde{e(1)}} &&{\widetilde{e(2)}} &&{\widetilde{e(3)}} &\cdots
&{\widetilde{e(r)}}
}\endxy$$
}

\noindent where the $e(i)$ and $\widetilde{e(i)}$ denote vertices and, for
each $i$, the
$p_i$, $q_i$ denote paths of positive length starting in distinct arrows.
That a module $M\in\lamod$ has the shown graph relative to a sequence
$x_1,\dots,x_r$ of top elements in particular encodes the following
information:
$x_i= e(i)x_i$, the $x_i$ are $K$-linearly independent modulo
$JM$, each $p_i$ has starting point $e(i)$ and endpoint $\widetilde{e(i)}$,
each $q_i$ has starting point $e(i)$ and endpoint $\widetilde{e(i+1)}$, and
the multiples $q_ix_i=p_{i+1}x_{i+1}$, $1\le i\le r-1$, are $K$-linearly
independent elements of the socle of $M$. The information encoded in the
graph guarantees, moreover, that there are no `other' nonzero multiples of
the $x_i$ apart from those shown; more precisely, the only paths in $\kgam$
not annihilating the element $x_i\in M$ are the right subpaths of $p_i$ and
$q_i$.

\head 2. Covariant finiteness of $\pinf(\lamod)$\endhead

The main goal of this section will be the development and application of 
conditions which guarantee that a simple left $\Lambda$-module fails to have
a left $\pinf(\lamod)$-ap\-prox\-i\-ma\-tion. We will start by
recalling a result of Auslander and Reiten, including a short alternate
proof akin to the arguments of the introduction.

\proclaim{Proposition 1} \cite{1, Proposition 4.2} If $\lfindim \Lambda \le
1$, then
$\pinf(\lamod)$ is covariantly finite in $\lamod$. \endproclaim

\demo{Proof} Suppose $\lfindim \Lambda\le 1$. By \cite{4, Theorem 4.2},
it suffices to prove that each direct product of objects in $\pinf(\lamod)$
belongs to $\arrpinf(\lamod)$. Clearly, each direct product $M$ of
objects from $\pinf(\lamod)$ has projective dimension $\le1$ in $\laMod$,
and hence \cite{9, Observation 5} shows that $M$ is a direct limit of
finitely generated modules of finite projective dimension as required.
\qed\enddemo

We will see later in this section that the conclusion of Proposition 1
breaks down for algebras of finitistic dimension 2.

\definition{Example 2} \cite{11} This is the example exhibited by Igusa,
Smal\o{} and Todor\-ov to show that $\pinf(\lamod)$ may fail to be
contravariantly finite, even in the case of a special biserial algebra
$\Lambda$.

Let $\Gamma$ be the quiver

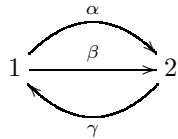
\ignore{
$$\xymatrixcolsep{4pc}
\xy\xymatrix{
1 \ar[r]^\beta \ar@/^1.5pc/[r]^\alpha &2 \ar@/^1.5pc/[l]^\gamma
}\endxy$$
}

\noindent and $I\subseteq \kgam$ such that, for $\Lambda= \kgam/I$, the
indecomposable projective left $\Lambda$-modules have the following graphs:

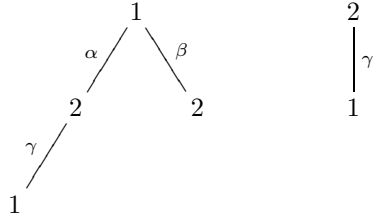
\ignore{
$$\xymatrixcolsep{1pc}
\xy\xymatrix{
 &&1 \edge[dl]_\alpha \edge[dr]^\beta &&&&2 \edge[d]^\gamma\\
 &2 \edge[dl]_\gamma &&2 &&&1\\
1}\endxy$$
}

\noindent Here $\lFindim \Lambda =1$, and hence
$\pinf(\lamod)$ is covariantly finite in $\lamod$ by Proposition 1.
\qed\enddefinition

In the sequel, $\A$ will denote a full subcategory of $\lamod$.

\proclaim{Criterion 3} Let $e(1),\dots,e(r)$ be vertices of $\Gamma$, and
let $p_1,\dots,p_r$, $q_1,\dots, \allowmathbreak q_r$ be $2r$ paths of
positive length in
$\kgam$, none of which is a subpath of any of the others. Moreover, suppose
that the following conditions {\rm (1)} and {\rm (2)} are satisfied:

{\rm (1)} For each natural number $n$, there exists a module $M_n\in\A$
having graph

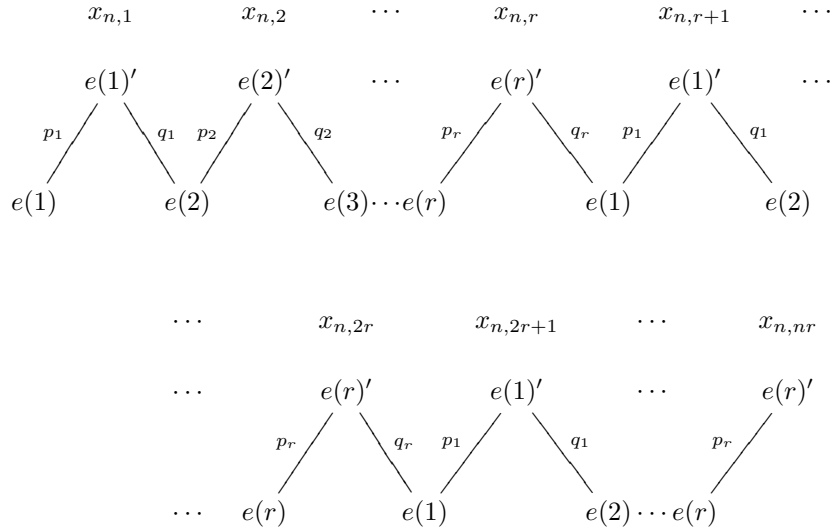
\ignore{
$$\xymatrixcolsep{0.35pc}\xymatrixrowsep{1pc}
\xy\xymatrix{
 &x_{n,1} &&x_{n,2} &\ar@{}[r]|{\displaystyle{\cdots}} &&x_{n,r} &&x_{n,r+1}
&\ar@{}[r]|{\displaystyle{\cdots}}&\\
 &e(1)' \edge[ddl]_{p_1} \edge[ddr]^{q_1} & &e(2)' \edge[ddl]_{p_2}
\edge[ddr]^{q_2} &\ar@{}[r]|{\displaystyle{\cdots}} &&e(r)' \edge[ddl]_{p_r}
\edge[ddr]^{q_r} &&e(1)' \edge[ddl]_{p_1} \edge[ddr]^{q_1}
&\ar@{}[r]|{\displaystyle{\cdots}}&\\
\\
e(1) &&e(2) &&e(3) \ar@{}[r]|{\displaystyle{\cdots}} &e(r) &&e(1) &&e(2)\\
\\
&&\cdots &&x_{n,2r} &&x_{n,2r+1} &\ar@{}[r]|{\displaystyle{\cdots}}
&&x_{n,nr}\\ 
&&\cdots &&e(r)' \edge[ddl]_{p_r} \edge[ddr]^{q_r} &&e(1)' \edge[ddl]_{p_1}
\edge[ddr]^{q_1} &\ar@{}[r]|{\displaystyle{\cdots}} &&e(r)' \edge[ddl]_{p_r}\\
\\
&&\cdots &e(r) &&e(1) &&e(2) \ar@{}[r]|{\displaystyle{\cdots}} &e(r)
}\endxy$$
}

\noindent relative to a suitable sequence of top elements $x_{n,1}, \dots,
x_{n,nr}$ of $M_n$ which are $K$-linearly independent modulo $JM_n$.

{\rm (2)} Each module $A\in\A$ has the following properties:

\qquad {\rm (i)} $e(1)(\Soc A) \subseteq p_1A$;

\qquad {\rm (ii)} $q_iA\cap (\Soc A)\subseteq p_{i+1}A$ for $i<r$, and
$q_rA\cap (\Soc A) \subseteq p_1A$;

\qquad {\rm (iii)} If $x\in A$ with $p_ix\in \Soc A$, then $q_ix\in \Soc A$.

Then $S= \Lambda e(1)/Je(1)$ fails to have a left $\A$-ap\-prox\-i\-ma\-tion.
\endproclaim

\demo{Proof} For $i>r$, let $s(i)$ be the integer in $\{1,\dots,r\}$ with
$i\equiv s(i) \pmod{r}$, and define $p_i := p_{s(i)}$, $q_i := q_{s(i)}$.

We assume that, to the contrary of our claim,  $S= \Lambda e(1)/Je(1)$ does
have a left $\A$-ap\-prox\-i\-ma\-tion $f : S\rightarrow A$ with $A\in\A$. Choose $n
>\dim_K A$ and write $x_1,\dots,x_{nr}$ for the elements $x_{n,1}, \dots,
x_{n,nr}$ of $M_n$. Moreover, let $g : S\rightarrow M_n$ be the embedding
which sends $e(1)+Je(1)$ to $p_1x_1$, and choose a homomorphism $h :
A\rightarrow M_n$ with $g= h\circ f$. Since $f(e(1)+Je(1)) \in e(1)\Soc A$,
condition 2(i) permits us to pick an element $a_1\in A$ such that
$f(e(1)+Je(1)) =p_1a_1$. In view of the equality $hf(e(1)+Je(1)) =p_1x_1$,
we see that $h(a_1)= x_1+y_1$ with $y_1\in \sum_{j \not\equiv 1 \pmod{r}}
\Lambda x_j +\sum_{j\equiv 1 \pmod{r}} Jx_j$.  Keep in mind that $p_1$
equals $p_{r+1}= \dots= p_{lr+1}$, but is not a subpath of the other $p_i$
or any of the $q_i$. Since $q_1$ is not a subpath of any of the paths
$p_1,\dots,p_r$, $q_2,\dots,q_r$ either, we infer that $q_1h(a_1)= q_1x_1$.
Due to our choice of $a_1$ such that $p_1a_1\in \Soc A$, condition 2(iii)
guarantees that $q_1a_1\in \Soc A$ as well, in other words, $q_1a_1\in
q_1A\cap \Soc A$. Next, condition 2(ii) yields $a_2\in A$ with $p_2a_2=
q_1a_1$. In view of $h(p_2a_2)= q_1h(a_1)= q_1x_1$, we obtain $h(a_2)=
x_2+y_2$ with $y_2\in \sum_{j\not\equiv 2 \pmod{r}} \Lambda x_j
+\sum_{j\equiv 2 \pmod {r}} Jx_j$; indeed this follows at once from the
nature of the graph of $M_n$ and the hypothesis that $p_2$ is not a subpath
of any of $p_1,p_3,\dots, p_r$, $q_1,\dots,q_r$. Now the non-occurrence of
$q_2$ as a subpath of $p_1,\dots,p_r$, $q_1,q_3, \dots,q_r$ allows us to
deduce $q_2h(a_2)= q_2x_2$, and since $p_2a_2= q_1a_1\in \Soc A$, we observe
that also $q_2a_2\in \Soc A$ by 2(iii). Thus 2(ii) in turn provides us with
an element $a_3\in A$ such that $p_3a_3= q_2a_2$. As above, we argue that
$q_3h(a_3)= q_3x_3$, and proceeding inductively, we thus obtain a sequence
$a_1,\dots,a_{nr}$ of elements in $A$ with the property that $q_ih(a_i)=
q_ix_i$ for $1\le i\le nr$. Since the elements $q_1x_1, \dots, q_{nr}x_{nr}$
of $M_n$ are $K$-linearly independent by hypothesis, so are
$a_1,\dots,a_{nr}$. But this contradicts our choice of $n$ exceeding
$\dim_KA$ and proves the criterion. \qed\enddemo

As our argument makes clear, the requirement that none of the $p_i$, $q_i$
is a subpath of any other as called for in the criterion can certainly be
weakened. Other than that, the chain reaction exhibited in the proof of
Criterion 3 appears prototypical for the failure of a simple module
$S$ to have a left
$\pinf(\lamod)$-ap\-prox\-i\-ma\-tion. In fact, specializing to special biserial
algebras, we believe the answer to the following question to be positive.

\definition{Problem 4} Suppose $S= \Lambda e(1)Je(1)$ is a simple left
module over a finite dimensional special biserial algebra $\Lambda=
\kgam/I$. If $S$ does not have a left $\pinf(\lamod)$-ap\-prox\-i\-ma\-tion, do
there exist paths $p_1,\dots,p_r$, $q_1,\dots,q_r$ in $\kgam$ satisfying the
conditions (1) and (2) of Criterion 3? \enddefinition

On the other hand, we point out that it is not known whether failure of
covariant finiteness of $\pinf(\lamod)$ in $\lamod$ implies that one of the
\underbar{simple} left $\Lambda$-modules is devoid of a left
$\pinf(\lamod)$-ap\-prox\-i\-ma\-tion. This is in contrast with our
level of information on contravariant finiteness of $\pinf(\lamod)$. Indeed,
as was shown by Auslander and Reiten in \cite{1},
$\pinf(\lamod)$ is contravariantly finite in $\lamod$ provided that all simple
$\Lambda$-modules have right $\pinf(\lamod)$-ap\-prox\-i\-ma\-tions. We therefore
propose

\definition{Problem 5} Is $\pinf(\lamod)$ covariantly finite in $\lamod$ in
case each simple left $\Lambda$-module has a left
$\pinf(\lamod)$-ap\-prox\-i\-ma\-tion? \enddefinition

Criterion 3 is particularly suited to the situation where $\A=
\pinf(\lamod)$ and $\Lambda$ is a finite dimensional special biserial
algebra (see \S1 for a definition). The representation theory of the
finitely generated modules over such algebras is exceptionally well
understood (see \cite{3,5,6,12,14}). In fact, the indecomposable objects
of
$\lamod$ fall into two classes, `bands' and `strings'. All we presently need
to know about bands is the following: If $B\in \lamod$ is a band, $e\in
\Lambda$ a primitive idempotent and $x\in e(\Soc B)$, then 
$x\in \sum_i u_iB\cap v_iB$ for pairs $(u_i,v_i)$ of paths of positive length
ending in $e$ such that, moreover, $u_i= \alpha_iu_i'$ and $v_i=
\beta_iv_i'$ with distinct arrows $\alpha_i$ and $\beta_i$. The strings, on
the other hand, are precisely the objects in $\lamod$ having graphs of the
form

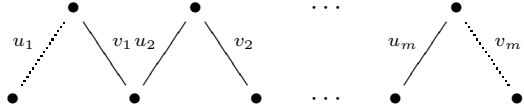
\ignore{
$$\xymatrixcolsep{1pc}
\xy\xymatrix{
 &\bullet \ar@{.}[dl]_{u_1} \edge[dr]^{v_1} &&\bullet \edge[dl]_{u_2}
\edge[dr]^{v_2} &&\cdots &&\bullet \edge[dl]_{u_m} \ar@{.}[dr]^{v_m}\\
\bullet &&\bullet &&\bullet &\cdots &\bullet &&\bullet
}\endxy$$
}

\noindent relative to a suitable sequence of top elements, such that, for
each $i$, the paths $u_i$, $v_i$ start in distinct arrows, and the partners
of each pair $(v_i,u_{i+1})$ end in distinct arrows.

\definition{Example 6} A finite dimensional special biserial algebra
$\Lambda$ for which $\pinf(\lamod)$ fails to be covariantly finite in
$\lamod$: Let $\Lambda= \kgam/I$, where $\Gamma$ is the quiver

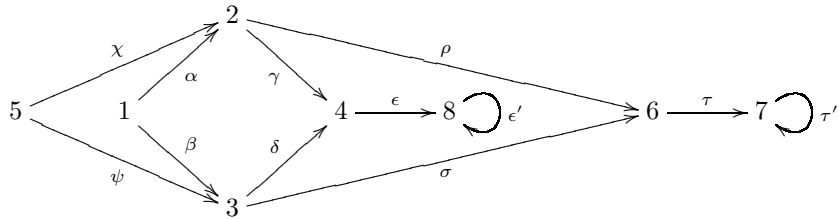
\ignore{
$$\xymatrixcolsep{2.5pc}
\xy\xymatrix{
 &&2 \ar[dr]_\gamma \ar[drrrr]^\rho\\
5 \ar[urr]^\chi \ar[drr]_\psi &1 \ar[ur]_\alpha \ar[dr]^\beta
&&4 \ar[r]^\epsilon &8 \ar@(ur,dr)^{\epsilon'} &&6 \ar[r]^\tau &7
\ar@(ur,dr)^{\tau'}\\
 &&3 \ar[ur]^\delta \ar[urrrr]_\sigma
}\endxy$$
}

\noindent and $I$ is generated by 
$$\gamma\alpha -\delta\beta,\ \rho\chi -\sigma\psi,\ \gamma\chi,\
\rho\alpha,\ \delta\psi,\ \sigma\beta,\ \epsilon\gamma,\ \epsilon'\epsilon,\
(\epsilon')^2,\ \tau'\tau,\ (\tau')^2,\ \tau\sigma.$$
Then the graphs of the indecomposable projective left $\Lambda$-modules are

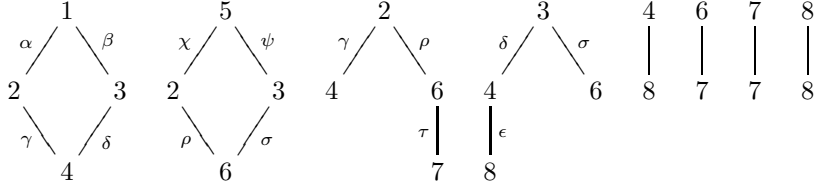
\ignore{
$$\xymatrixcolsep{0.75pc}\xymatrixrowsep{1.5pc}
\xy\xymatrix{
 &1 \edge[dl]_\alpha \edge[dr]^\beta &&&5 \edge[dl]_\chi \edge[dr]^\psi
&&&2 \edge[dl]_\gamma \edge[dr]^\rho &&&3 \edge[dl]_\delta
\edge[dr]^\sigma &&4 \edge[d] &6 \edge[d] &7 \edge[d] &8 \edge[d]\\
2 \edge[dr]_\gamma &&3 \edge[dl]^\delta &2 \edge[dr]_\rho &&3
\edge[dl]^\sigma &4 &&6 \edge[d]_\tau &4 \edge[d]^\epsilon &&6 &8 &7 &7
&8\\
 &4 &&&6 &&&&7 &8
}\endxy$$
}

\noindent We will use Criterion 3 to show that $S= \Lambda e_3/Je_3$ does
not have a left $\pinf(\lamod)$-ap\-prox\-i\-ma\-tion. For that purpose, let $r=2$,
$e(1)= e_3$, $e(2)= e_5$, $p_1=\beta$, and $p_2=\chi$, $q_1=\alpha$,
$q_2=\psi$. Then condition 1 of our criterion is satisfied by the modules
$M_n$ with graphs

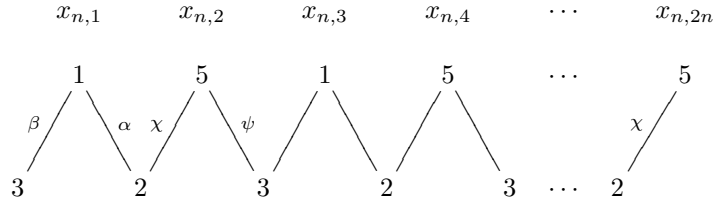
\ignore{
$$\xymatrixcolsep{0.5pc}\xymatrixrowsep{1pc}
\xy\xymatrix{
 &x_{n,1} &&x_{n,2} &&x_{n,3} &&x_{n,4} &&\cdots &&x_{n,2n}\\
 &1 \edge[ddl]_\beta \edge[ddr]^\alpha &&5 \edge[ddl]_\chi \edge[ddr]^\psi &&1
\edge[ddl] \edge[ddr] &&5 \edge[ddl] \edge[ddr] &&\cdots &&5 \edge[ddl]_\chi\\
\\
3 &&2 &&3 &&2 &&3 &\cdots &2
}\endxy$$
}

\noindent Indeed, one readily verifies that the first syzygy $\Omega^1(M_n)$
of $M_n$ has graph

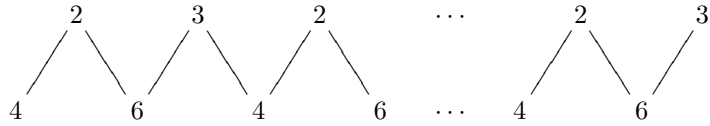
\ignore{
$$\xymatrixcolsep{1pc}
\xy\xymatrix{
 &2 \edge[dl] \edge[dr] &&3 \edge[dl] \edge[dr] &&2 \edge[dl] \edge[dr]
&&\cdots &&2 \edge[dl] \edge[dr] &&3 \edge[dl]\\
4 &&6 &&4 &&6 &\cdots &4 &&6
}\endxy$$
}

\noindent and the second syzygy $\Omega^2(M_n)$ is a direct sum of modules
with graphs

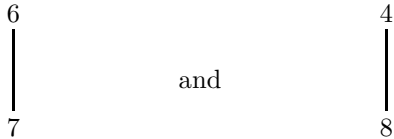
\ignore{
$$\xymatrixrowsep{0.75pc}
\xy\xymatrix{
6 \edge[dd] &&&&4 \edge[dd]\\
 &&\txt{and}\\
7 &&&&8
}\endxy$$
}

\noindent This shows that $\pdim M_n =2$; in particular, $M_n\in
\pinf(\lamod)$.

It is obvious that for every module $A\in\lamod$ and every element $x\in A$
the implications ($\beta x\in \Soc A$ $\implies$ $\alpha x\in \Soc A$) and
$(\chi x\in \Soc A$
$\implies$ $\psi x\in \Soc A$) hold. So condition 2(iii) of Criterion 3 is
met as well.

To check the remaining conditions under (2), we start by observing that the
only nontrivial paths of $\kgam$ ending in $e_3$ are the arrows $\beta$ and
$\psi$, while the only nontrivial paths ending in $e_2$ are the arrows
$\alpha$ and $\chi$. Given our comments on bands, this implies that, for
each band $B\in \lamod$, we have $e_3(\Soc B) \subseteq \beta B\cap \psi B$
and, in particular, $\psi B\cap (\Soc B)\subseteq \beta B$. Analogously,
$\alpha B\cap (\Soc B) \subseteq e_2(\Soc B) \subseteq \chi B$. Hence, we
may focus our attention on the situation where $A\in \pinf(\lamod)$ is a
string with $S_3\subseteq \Soc A$. Note that $A$ cannot be simple, since
$\pdim S_3 =\infty$. The only possibility for a copy of $S_3 \subseteq \Soc
A$ not to belong to $\beta A$ is that of a string $A$ having graph

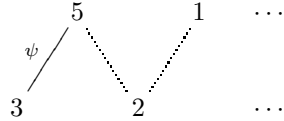
\ignore{
$$\xymatrixcolsep{1pc}
\xy\xymatrix{
 &5 \edge[dl]_\psi \ar@{.}[dr] &&1 \ar@{.}[dl] &\cdots\\
3 &&2 &&\cdots
}\endxy$$
}

\noindent In that case $\Omega^1(A)$ would have a graph of one of the
following types

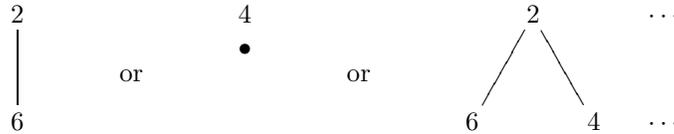
\ignore{
$$\xymatrixcolsep{1pc}\xymatrixrowsep{0.75pc}
\xy\xymatrix{
2 \edge[dd] &&&&4 \save+<0ex,-3ex> \drop{\bullet} \restore &&&&&2 \edge[ddl]
\edge[ddr] &&\cdots\\
 &&\txt{or} &&&&\txt{or}\\
6 &&&& &&&&6 &&4 &\cdots
}\endxy$$
}

\noindent and $\Omega^2(A)$ would be a copy of $S_4\oplus S_7$ in the first
case, a copy of $S_8$ in the second, and have a direct summand isomorphic to
$S_7$ in the third. In all of these cases, we would have $\pdim \Omega^2(A)=
\infty$ contradicting our choice of $A$ in $\pinf(\lamod)$. This proves
$e_3(\Soc A) \subseteq \beta A$ and thus 2(i).

For 2(ii) it suffices to observe that any $A\in \pinf(\lamod)$ has the
stronger property that $e_2(\Soc A) \subseteq \chi A$.
The argument is analogous to the one we just completed for $e_3$ and $\beta$.

Thus, by Criterion 3, $S= \Lambda e_3/Je_3$ does not have a left
$\pinf(\lamod)$-ap\-prox\-i\-ma\-tion. \qed\enddefinition

We remark that $\lFindim \Lambda =2$ in the preceding example, which shows
that Auslander-Reiten's Proposition 1 cannot be extended to the case of
finitistic dimensions exceeding 1.

\head 3. Curing failure of contravariant finiteness of
$\pinf(\lamod)$\endhead

We conjecture that, for any monomial relation algebra and for any special
biserial algebra $\Lambda$, there exists a sequence of one-point extensions
$\Lambda=
\Lambda_0, \dots, \Lambda_m= \Delta$ such that $\pinf(\dlmod)$ is
contravariantly finite in $\dlmod$. At this point, our conviction is based
mainly on a long list of examples. One of our most interesting examples
shows that, in general, one cannot expect $\Delta$ to retain the `good'
properties of $\Lambda$ in this process, in other words, a `cure' for
failure of contravariant finiteness of $\pinf(\lamod)$ by successive
one-point extensions, will usually involve a trade-off. Here we will
construct an example of a monomial relation algebra $\Lambda$, together with
a sequence $\Lambda= \Lambda_0, \Lambda_1, \Lambda_2= \Delta$ as above such
that $\Delta$ cannot be chosen within the class of monomial relation
algebras. Simultaneously, this example will illustrate the potential
intricacy of the structure of the minimal right
$\pinf(\dlmod)$-ap\-prox\-i\-ma\-tions of the simple modules.

\definition{Example 7} Let $\Lambda= \kgam/I$, where $\Gamma$ is the
quiver

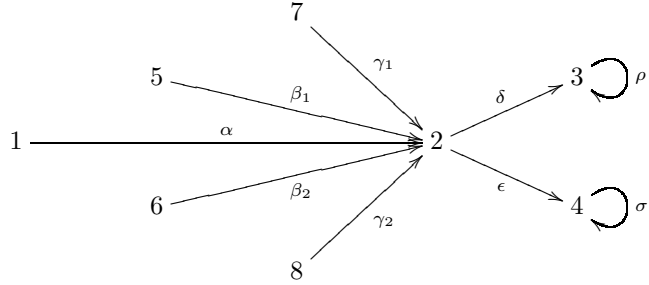
\ignorethree{
$$\xymatrixcolsep{3.5pc}\xymatrixrowsep{1pc}
\xy\xymatrix{
 &&7 \ar[ddr]^{\gamma_1}\\
 &5 \ar[drr]^{\beta_1} &&&3 \ar@(ur,dr)^\rho\\
1 \ar[rrr]^\alpha &&&2 \ar[ur]^\delta \ar[dr]_\epsilon\\
 &6 \ar[urr]_{\beta_2} &&&4 \ar@(ur,dr)^\sigma\\
 &&8 \ar[uur]_{\gamma_2}
}\endxy$$
}

\noindent and the ideal $I\subseteq \kgam$ is generated by
$$\delta\alpha,\ \epsilon\alpha,\ \epsilon\beta_i\ (i=1,2),\ \delta\gamma_i\
(i=1,2),\ \rho\delta,\ \sigma\epsilon,\ \rho^2,\ \sigma^2.$$
This yields indecomposable projective left $\Lambda$-modules with graphs

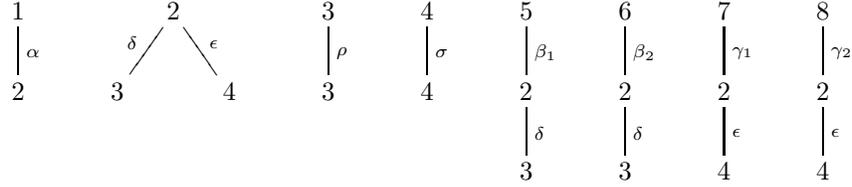
\ignore{
$$\xymatrixcolsep{0.85pc}\xymatrixrowsep{1.5pc}
\xy\xymatrix{
1 \edge[d]^\alpha &&&2 \edge[dl]_\delta \edge[dr]^\epsilon &&&3
\edge[d]^\rho &&4 \edge[d]^\sigma &&5 \edge[d]^{\beta_1} &&6
\edge[d]^{\beta_2} &&7 \edge[d]^{\gamma_1} &&8 \edge[d]^{\gamma_2}\\
2 &&3 &&4 &&3 &&4 &&2 \edge[d]^\delta &&2 \edge[d]^\delta &&2
\edge[d]^\epsilon &&2 \edge[d]^\epsilon\\
 &&&&&&&&&&3 &&3 &&4 &&4
}\endxy$$
}

\noindent To see that $\pinf(\lamod)$ fails to be contravariantly finite in
$\lamod$, more precisely, that $\Lambda e_1/Je_1$ fails to have a right
$\pinf(\lamod)$-ap\-prox\-i\-ma\-tion, we exhibit a $\pinf(\lamod)$-phan\-tom
of infinite $K$-dimension for $\Lambda e_1/Je_1$ (see \cite{7, Definition 5
and Theorem 9}). The routine check that the following module $H$ is indeed
such a phantom is left to the reader: $H= \varinjlim H_n$ where
$$H= \biggl( \Lambda z \oplus \bigoplus_{i=1}^n \Lambda x_i \oplus
\bigoplus_{i=1}^n \Lambda y_i \biggr) \biggm/ U_n,$$
with $z=e_1$, $x_{2m-1}= e_5$, $x_{2m}= e_6$, $y_{2m-1}= e_7$, $y_{2m}= e_8$
for
$m\ge1$, and
$$U_n= \Lambda (\alpha z-\beta_1x_1 -\gamma_1y_1) +\sum_{i=1}^{n-1} \Lambda
(\tilgam_iy_i -\tilbet_{i+1}x_{i+1} -\tilgam_{i+1}y_{i+1}),$$
with $\tilgam_i$ equal to $\gamma_1$ or $\gamma_2$, depending on whether $i$
is odd or even, and $\tilbet_i$ equal to $\beta_1$ or $\beta_2$, depending
on whether $i$ is odd or even. The modules $H_n$ can be pictured via graphs
of the form

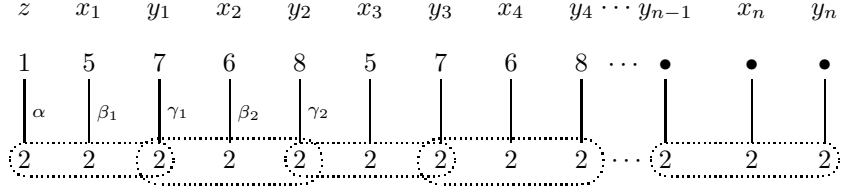
\ignore{
$$\xymatrixcolsep{0.9pc}\xymatrixrowsep{0.75pc}
\xy\xymatrix{
z &x_1 &y_1 &x_2 &y_2 &x_3 &y_3 &x_4 &y_4
\ar@{}[r]|{\displaystyle{\cdots}\phantom{\cdot\cdot}} &y_{n-1} &x_n &y_n\\ 
1 \edge[dd]^\alpha &5 \edge[dd]^{\beta_1} &7 \edge[dd]^{\gamma_1} &6
\edge[dd]^{\beta_2} &8 \edge[dd]^{\gamma_2} &5 \edge[dd] &7 \edge[dd] &6
\edge[dd] &8 \edge[dd] \ar@{}[r]|{\displaystyle{\cdots}} &\bullet \edge[dd]
&\bullet \edge[dd] &\bullet \edge[dd]\\
\\
2 \save[0,0];[0,2]**\frm<6pt>{.}\restore &2 &2
\save[0,0]+(-3,-3);[0,2]+(3,3)**\frm<10pt>{.}\restore &2 &2
\save[0,0];[0,2]**\frm<6pt>{.}\restore &2 &2
\save[0,0]+(-3,-3);[0,2]+(3,3)**\frm<10pt>{.}\restore &2 &2
\ar@{}[r]|{\phantom{\cdot}\displaystyle{\cdots}} &2
\save[0,0];[0,2]**\frm<6pt>{.}\restore &2 &2  
}\endxy$$
}

\noindent relative to the top elements $z,x_1,\dots,x_n,y_1,\dots,y_n$,
where we extend our graphing conventions as follows: The dotted loop around
the vertices labeled `2' which represent $\alpha z$,
$\beta_1x_1$, $\gamma_1y_1$ indicates that any two of the three listed
vectors are
$K$-linearly independent while
$\dim_K (K\alpha z+K\beta_1x_1 +K\gamma_1y_1) =2$. The same holds for the
additional triples
$\tilgam_iy_i, \tilbet_ix_{i+1}, \tilgam_{i+1}y_{i+1}$ surrounded by loops.
Note that $\Omega^1(H_n)\cong (\Lambda e_2)^n$ for $n\in\NN$; in particular,
$H_n\in \pinf(\lamod)$ for $n\in\NN$.

Let $\Gamma_1$ be the quiver obtained from $\Gamma$ by adding a single
vertex, labeled 9, and two arrows leaving 9, namely $\chi_1 : 9\rightarrow
5$ and $\chi_2 : 9\rightarrow 6$. Moreover, let $\Lambda_1= \kgam_1/I_1$,
where the ideal $I_1\subseteq \kgam_1$ is generated by $I$ and the relation
$\beta_1\chi_1 -\beta_2\chi_2$.

Next let $\Gamma_2$ be the quiver obtained from $\Gamma_1$ by adding a
single vertex, 10, and two arrows leaving 10, namely $\psi_1 : 10\rightarrow
7$ and $\psi_2 : 10\rightarrow 8$. Now $\Delta= \Lambda_2= \kgam_2/I_2$,
where $I_2$ is generated by $I_1$ and the relation $\gamma_1\psi_1
-\gamma_2\psi_2\in \kgam_2$.

Clearly, $\Lambda_1$ is a one-point extension of $\Lambda= \Lambda_0$, and
$\Delta= \Lambda_2$ is a one-point extension of $\Lambda_1$. The `new'
indecomposable projective left $\Delta$-modules have graphs

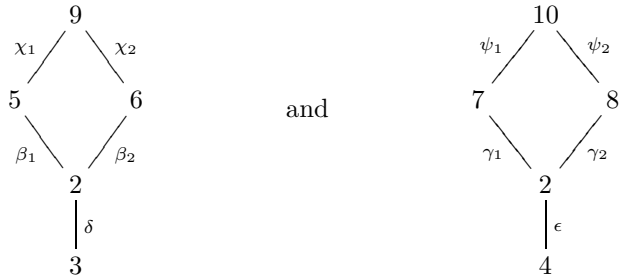
\ignorethree{
$$\xymatrixcolsep{1pc}\xymatrixrowsep{1.5pc}
\xy\xymatrix{
 &9 \edge[dl]_{\chi_1} \edge[dr]^{\chi_2} &&&& &&&&10 \edge[dl]_{\psi_1}
\edge[dr]^{\psi_2}\\
5 \edge[dr]_{\beta_1} &&6 \edge[dl]^{\beta_2} &&&\txt{and} &&&7
\edge[dr]_{\gamma_1} &&8 \edge[dl]^{\gamma_2}\\
 &2 \edge[d]^\delta &&&& &&&&2 \edge[d]^\epsilon\\
 &3 &&&& &&&&4
}\endxy$$
}

\noindent Note that $\Delta e_i= \Lambda e_i$ for $1\le i\le 8$.

One can verify that $\pinf(\dlmod)$ is contravariantly finite in $\dlmod$ by
exhibiting right $\pinf(\dlmod)$-ap\-prox\-i\-ma\-tions of the simple left
$\Delta$-modules $S_i= \Delta e_i/J(\Delta)e_i$ for $1\le i\le 10$. It is
comparatively easy to see that the following are the (minimal) right
$\pinf(\dlmod)$-ap\-prox\-i\-ma\-tions of $S_2,\dots, S_{10}$: Namely,
$\Delta e_i$ for $i= 2,3,4$, and 

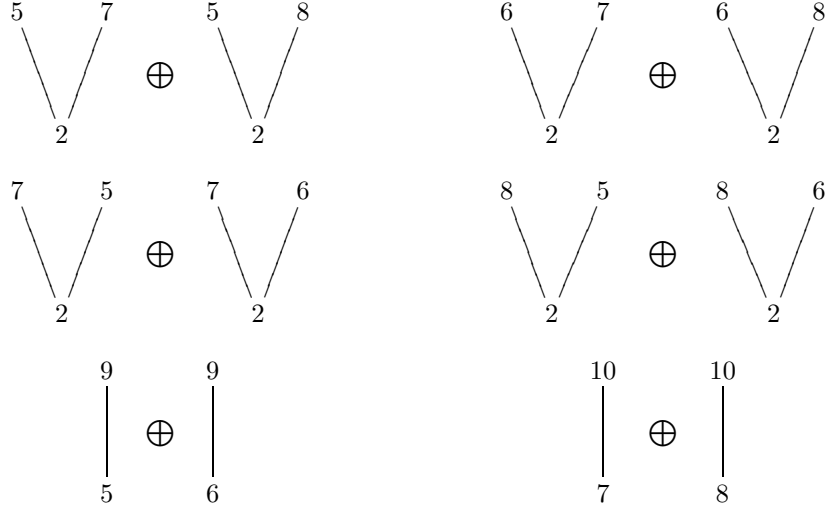
\ignorethree{
$$\xymatrixcolsep{0.5pc}\xymatrixrowsep{0.75pc}
\xy\xymatrix{
5 \edge[ddr] &&7 \edge[ddl] &&5 \edge[ddr] &&8 \edge[ddl] &&&&&&6 \edge[ddr]
&&7 \edge[ddl] &&6 \edge[ddr] &&8 \edge[ddl]\\
 &&&\bigoplus &&&&&& &&&&&&\bigoplus\\
 &2 &&&&2 &&&&&&&&2 &&&&2\\
7 \edge[ddr] &&5 \edge[ddl] &&7 \edge[ddr] &&6 \edge[ddl] &&&&&&8 \edge[ddr]
&&5 \edge[ddl] &&8 \edge[ddr] &&6 \edge[ddl]\\
 &&&\bigoplus &&&&&& &&&&&&\bigoplus\\
 &2 &&&&2 &&&&&&&&2 &&&&2\\
 &&9 \edge[dd] &&9 \edge[dd] &&&&&&&&&&10 \edge[dd] &&10 \edge[dd]\\
 && &\bigoplus &&&&&& &&&&&&\bigoplus\\
 &&5 &&6 &&&&&&&&&&7 &&8
}\endxy$$
}

\noindent for $i=5,6,\dots,10$, respectively. We will sketch an argument
backing the claim that the canonical epimorphism $f : A_1\rightarrow S_1$
with 
$$A_1= (\Delta e_1\oplus \Delta e_9\oplus \Delta e_{10}) /\Delta( \alpha,
\beta_1\chi_1, \gamma_1\psi_1)$$
is a right $\pinf(\dlmod)$-ap\-prox\-i\-ma\-tion of $S_1$. To buttress
intuition, start by noting that $A_1$ has graph

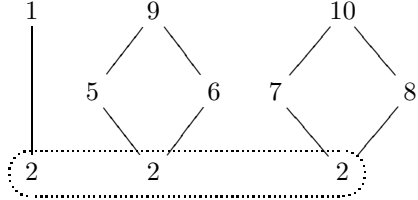
\ignorethree{
$$\xymatrixcolsep{1pc}\xymatrixrowsep{1.5pc}
\xy\xymatrix{
1 \edge[dd] &&9 \edge[dl] \edge[dr] &&&10 \edge[dl] \edge[dr]\\
 &5 \edge[dr] &&6 \edge[dl] &7 \edge[dr] &&8 \edge[dl]\\
2 \save[0,0]+(-3,-3);[0,5]+(3,3)**\frm<10pt>{.}\restore &&2 &&&2
}\endxy$$
}

\noindent with the above convention for loops. Observe, moreover, that
$\Omega^1_\Delta(A_1)\cong \Delta e_2$, whence $A_1\in \pinf(\dlmod)$.

We will sketch an argument showing that each epimorphism $g : M\rightarrow
S_1$ with
$M\in
\pinf(\dlmod)$ factors through $f$. It is clearly harmless to assume that
$M$ is indecomposable. Moreover, it suffices to consider the case where
$\ker(g)$ {\it does not contain any nonzero submodules in}
$\pinf(\dlmod)$; indeed, given $U\subseteq M$ with $U\in \pinf(\dlmod)$, it
is enough to factor the induced map ${\overline g} : M/U \rightarrow S_1$
through $f$. In particular, this means that $M$ does not contain any
submodules isomorphic to $\Delta e_i$, with $i\ge 2$, nor any submodules with
graphs of type

\ignorethree{
$$\xymatrixcolsep{1pc}\xymatrixrowsep{0.75pc}
\xy\xymatrix{
10 \edge[dd] &&10 \edge[dd] &&11 \edge[dd] &&11 \edge[dd]\\
 &\txt{or} &&\txt{or} &&\txt{or}\\
6 &&7 &&8 &&9
}\endxy$$
}

\noindent As a consequence, we can zero in on the structure of $M$ as
follows: Let $\A$ be the class of $\Delta$-modules isomorphic to $\Delta
e_1$, $\B$ the class of modules isomorphic to one of the $\Delta$-modules
with graphs

\ignorethree{
$$\xymatrixcolsep{1pc}\xymatrixrowsep{1pc}
\xy\xymatrix{
5 \edge[dd] &&6 \edge[dd] &&&9 \edge[dl] \edge[dr]\\
 &\txt{or} &&\txt{or} &5 \edge[dr] &&6 \edge[dl]\\
2 &&2 &&&2
}\endxy$$
}

\noindent and, finally, $\C$ the class of those $\Delta$-modules which have
one of the graphs

\ignorethree{
$$\xymatrixcolsep{1pc}\xymatrixrowsep{1pc}
\xy\xymatrix{
7 \edge[dd] &&8 \edge[dd] &&&10 \edge[dl] \edge[dr]\\
 &\txt{or} &&\txt{or} &7 \edge[dr] &&8 \edge[dl]\\
2 &&2 &&&2
}\endxy$$
}

\noindent Our assumptions on $M$ guarantee that, up to isomorphism, $M=
X/Y\in \pinf(\dlmod)$, where 
$$X= \bigoplus_{1\le i\le m_1} A_i \oplus \bigoplus_{1\le i\le m_2} B_i
\oplus \bigoplus_{1\le i\le m_3} C_i$$
with $A_i\in\A$, $B_i\in\B$, $C_i\in\C$ and $Y\subseteq \Soc X\cong
S_2^{m_1+m_2+m_3}$.

Let $a_i$, $b_i$, $c_i$ be top elements of $A_i$, $B_i$, $C_i$,
respectively, and let $a'_i= \alpha a_i$, $b'_i= (\beta_1 +\beta_2
+\beta_1\chi_1)b_i$, $c'_i= (\gamma_1+\gamma_2 +\gamma_1\psi_1)c_i$.
Moreover, let
$$\bigoplus_{1\le i\le m_1} \Delta \abar_i \oplus \bigoplus_{1\le i\le m_2}
\Delta \bbar_i \oplus \bigoplus_{1\le i\le m_3} \Delta \cbar_i
\longrightarrow X/Y$$
be the obvious projective cover of $X/Y$ mapping $\abar_i$, $\bbar_i$,
$\cbar_i$ to $a_i$, $b_i$, $c_i$ respectively. Then
$$\Soc X= \bigoplus_{1\le i\le m_1} \Delta a'_i \oplus \bigoplus_{1\le i\le
m_2} \Delta b'_i \oplus \bigoplus_{1\le i\le m_3} \Delta c'_i,$$
and we can therefore write $Y$ in the form $\bigoplus_{1\le h\le t} \Delta
y_h$ with 
$$y_h= \sum_{1\le i\le m_1} k_{hi}a'_i +\sum_{1\le i\le m_2} l_{hi}b'_i
+\sum_{1\le i\le m_3} m_{hi} c'_i,$$
where $k_{hi}, l_{hi}, m_{hi} \in K$ such that $\Delta y_h\cong \Delta e_2=
\Lambda e_2$ for each $h$. Let $\bfk_h= (k_{hi})_{1\le i\le m_1} \in
K^{m_1}$, $\bfl_h= (l_{hi})_{1\le i\le m_2} \in K^{m_2}$, and $\bfm_h=
(m_{hi})_{1\le i \le m_3}$ in $K^{m_3}$. Then the vectors $\bfl_1,\dots,
\bfl_t\in K^{m_2}$ are $K$-linearly independent; indeed if we had
$\sum_{1\le h\le t} d_h\bfl_h =0$ with $d_h\in K$ not all zero, we would
obtain a top element $z$ of $\Omega^1_\Delta(M)$
with the property that $\delta z=0$, namely $z= \sum_{1\le h\le t} d_hz_h$,
where 
$$\align z_h= \sum_{1\le i\le m_1} k_{hi}\alpha\abar_i &+\sum_{1\le i\le
m_2} l_{hi}(\beta_1 +\beta_2 +\beta_1\chi_1)\bbar_i \\
 &+\sum_{1\le i\le m_3}
m_{hi} (\gamma_1 +\gamma_2 +\gamma_1\psi_1)\cbar_i \endalign$$ 
in $\Omega^1_\Delta(M)$. This would place a direct
summand isomorphic to
$S_3$ into
$\Omega^1_\Delta(M)$, which -- in view of
$\pdim_\Delta \Delta e_3/J(\Delta e_3)= \pdim \Lambda e_3/Je_3 =\infty$ --
is incompatible with $\pdim_\Delta (M) <\infty$. Similarly $\bfm_1,\dots,
\bfm_t$ in $K^{m_3}$ are linearly independent, since otherwise we would
obtain a direct summand isomorphic to $S_4$ in $\Omega^1_\Delta(M)$. If we
set $\tlg(a_i)= \overline{(r_ie_1,0,0)} \in A_1$, where $f(a_i)=
\overline{r_ie_1}$ with $r_i\in K$, the above independence information allows
us to extend the assignment
$\tlg$ to a homomorphism
$\tlg : M\rightarrow A$. Any such homomorphism clearly satisfies $f\circ \tlg
=g$.
\bigskip

Now let $\Lambda= R_0, R_1, \dots, R_m=R$ be successive one-point extensions
such that $R$ is a monomial relation algebra. We leave the justification of
our claim that $\pinf(\rmod)$ fails to be contravariantly finite in $\rmod$
as an exercise, but provide hints of the underlying ideas.

1.) Due to the fact that $R$ is a monomial relation algebra obtained from
$\Lambda$ via one-point extensions, the following holds: Whenever a left
$R$-module
$M$ has a submodule
$N$ with graph

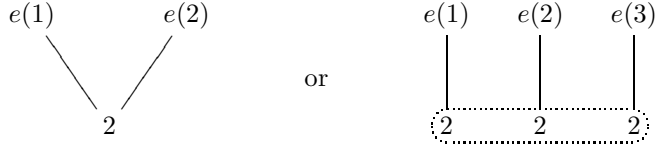
\ignorethree{
$$\xymatrixcolsep{1pc}\xymatrixrowsep{0.75pc}
\xy\xymatrix{
e(1) \edge[ddr] &&e(2) \edge[ddl] &&&&e(1) \edge[dd] &e(2) \edge[dd] &e(3)
\edge[dd]\\
 &&&&\txt{or}\\
 &2 &&&&&2 \save[0,0];[0,2]**\frm<6pt>{.}\restore &2 &2 
}\endxy$$
}

\noindent where $e(1)$, $e(2)$, $e(3)$ are distinct vertices in
$\{e_1,e_5,e_6,e_7,e_8\}$, the first syzygy $\Omega^1_R(M)$ of $M$ has a top
element $x$ of type $e_2$; moreover, if $e(1)$, $e(2)$, $e(3)$ belong to
$\{e_1,e_5,e_6\}$, we can choose $x$ such that $Rx$ has graph

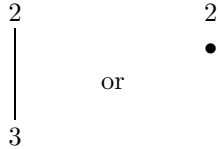
\ignorethree{
$$\xymatrixcolsep{2pc}\xymatrixrowsep{1pc}
\xy\xymatrix{
2 \edge[dd] &&2 \save+<0ex,-3ex> \drop{\bullet} \restore\\
 &\txt{or}\\
3
}\endxy$$
}

\noindent and if $e(1),e(2),e(3)\in \{e_1,e_7,e_8\}$, we can choose $x$ such
that $Rx$ has graph

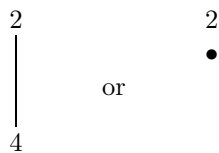
\ignorethree{
$$\xymatrixcolsep{2pc}\xymatrixrowsep{1pc}
\xy\xymatrix{
2 \edge[dd] &&2 \save+<0ex,-3ex> \drop{\bullet} \restore\\
 &\txt{or}\\
4
}\endxy$$
}

\noindent In each of these two situations, $\pdim_R(M) =\infty$.

2.) If there were a $\pinf(\rmod)$-approximation $B_1$ of $Re_1/J(R)e_1=
\Lambda e_1/Je_1= S_1$, then all homomorphisms in $\Hom_R(H_n,S_1)=
\Hom_\Lambda(H_n,S_1)$ with $H_n$ as above would factor through $B_1$,
because $H_n\in \pinf(\rmod)$. Using the first part, one would deduce the
existence of a submodule of $B_1$ with graph

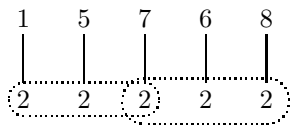
\ignorethree{
$$\xymatrixcolsep{1pc}\xymatrixrowsep{1.5pc}
\xy\xymatrix{
1 \edge[d] &5 \edge[d] &7 \edge[d] &6 \edge[d] &8 \edge[d]\\
2 \save[0,0];[0,2]**\frm<6pt>{.}\restore &2 &2
\save[0,0]+(-3,-3);[0,2]+(3,3)**\frm<10pt>{.}\restore &2 &2
}\endxy$$
}

\noindent and then proceed to show that $H= \varinjlim H_n$ would
still be a $\pinf(\rmod)$-phantom for $S_1$.

\enddefinition


\Refs

\ref\no1\by M. Auslander and I. Reiten \paper Applications of
contravariantly finite subcategories \jour Advances in Math. \vol 86 \yr
1991 \pages 111-152 \endref

\ref\no2\by M. Auslander and S.O. Smal\o \paper Preprojective modules over
Artin algebras \jour J. Algebra \vol 66 \yr 1980 \pages 61-122 \endref

\ref\no3\by M.C.R. Butler and C.M. Ringel \paper Auslander-Reiten sequences
with few middle terms and applications to string algebras \jour Comm.
Algebra \vol 15 \yr 1987 \pages 145-179 \endref

\ref\no4\by W. Crawley-Boevey \paper Locally finitely presented additive
categories \jour Comm. Algebra \vol 22 \yr 1994 \pages 1641-1674 \endref

\ref\no5\by P.W. Donovan and M.R. Freislich \paper The indecomposable
representations of certain groups with dihedral Sylow subgroups \jour Math.
Ann. \vol 238 \yr 1978 \pages 207-216 \endref

\ref\no6\by I.M. Gelfand and V.A. Ponomarev \paper Indecomposable
representations of the Lorentz group \jour Uspehi Mat. Nauk \vol 23 \yr 1968
\pages 3-60 \transl English Transl. \jour Russian Math. Surveys \vol 23
\yr 1969 \pages 1-58 \endref

\ref\no7\by D. Happel and B. Huisgen-Zimmermann \paper Viewing finite
dimensional representations through infinite dimensional ones \paperinfo
manuscript \endref

\ref\no8\by B. Huisgen-Zimmermann \paper Predicting syzygies over monomial
relation algebras \linebreak \jour manu\-scrip\-ta math. \vol 70 \yr 1991
\pages 157-182
\endref

\ref\no9\bysame \paper Homological assets of positively graded
representations of finite dimensional algebras \inbook in Representations of
Algebras (Ottawa 1992) \eds V. Dlab and H. Lenzing \bookinfo
	Canad. Math. Soc. Conf. Proc. Series 14 \yr 1993 \pages 463-475\endref

\ref\no10\by B. Huisgen-Zimmermann and S.O. Smal\o \paper A homological
bridge between finite and infinite dimensional representations of algebras
\paperinfo manuscript \endref

\ref\no11\by K. Igusa, S.O. Smal\o, and G. Todorov \paper Finite
projectivity and contravariant finiteness \jour Proc. Amer. Math. Soc. \vol
109 \yr 1990 \pages 937-941 \endref

\ref\no12\by C.M. Ringel \paper The indecomposable representations of the
dihedral 2-groups \jour Math. Ann. \vol 214 \yr 1975 \pages 19-34 \endref

\ref\no13\bysame \book Tame Algebras and Integral Quadratic Forms \bookinfo
Lecture Notes in Math. 1099 \publaddr Berlin \yr 1984 \publ Springer-Verlag
\endref

\ref\no14\by B. Wald and J. Waschb\"usch \paper Tame biserial algebras \jour
J. Algebra \vol 95 \yr 1985 \pages 480-500 \endref

\endRefs
\enddocument